\documentstyle[12pt]{article}
\newcommand{\la}[1]{\mbox{\large $#1$}}
\newcommand{\La}[1]{\mbox{\Large $#1$}}
\newcommand{\LA}[1]{\mbox{\LARGE $#1$}}

\newcommand{\mb}[1]{\mbox{#1}}

\newcommand{\lnab}[1]{\la{\nabla}_{\!\!#1}}

\newcommand{\lnabo}[1]{\la{\nabla}^{\bot}_{\!\!#1}} 

\newcommand{\qed}{\mbox{~~~\boldmath $\Box$}} \newcommand{\R}[1]{I\!\!R^{#1}}
\newcommand{\Co}{l\!\!\! C}
\newcommand{\ra}{\rightarrow}
\newcommand{\non}{\nonumber}

\newcommand{\al}{\alpha}
\newcommand{\be}{\beta}
\newcommand{\ga}{\gamma}
\newcommand{\bal}{\bar{\alpha}}
\newcommand{\bbe}{\bar{\beta}}
\newcommand{\bga}{\bar{\gamma}}
\newcommand{\bmu}{\bar{\mu}}
\newcommand{\brho}{\bar{\rho}}
\newcommand{\tg}{\tilde{g}}
\newcommand{\gdf}[3]{g(\mbox{\large ${\nabla}$}_{\!\!{#1}} dF({#2}),JdF({#3}))}
\newcommand{\gf}[3]{g_{#1}{#2}{#3}}
\newcommand{\Jw}{J_{\omega}}
\newcommand{\kw}{{\cal K}_{\omega}}

\newcommand{\Fw}{F^{*}\omega}
\newtheorem{Lm}{Lemma}
\newtheorem{Pp}{Proposition}
\newtheorem{Th}{Theorem}
\newtheorem{Cr}{Corollary}

\begin{document}
\baselineskip .6cm
\setcounter{page}{1}
\title{Broadly-Pluriminimal Submanifolds of K\"{a}hler-Einstein Manifolds} \author{
Isabel M.\ C.\ Salavessa\raisebox{.8ex}{\scriptsize 1} and 
Giorgio Valli\raisebox{.8ex}{\scriptsize \dag2} }
\protect\footnotetext{Deceased on October 2nd, 1999}
\date{}
\maketitle ~~~\\[-15mm]
{\footnotesize 1 Centro de F\'{\i}sica das
Interac\c{c}\~{o}es Fundamentais, Instituto Superior T\'{e}cnico,
Edif\'{\i}cio Ci\^{e}ncia,\\[-2mm] piso 3,
1096 LISBOA Codex, Portugal.~~
e-mail: isabel@cartan.ist.utl.pt}\\[-1mm] 
{\footnotesize 2 Dipartimento di
Matematica, Universit\`{a} di Pavia,
Via Abbiategrasso 215,
27100 PAVIA, Italy.}\\[1cm]
{\small {\bf Abstract:} We define broadly-pluriminimal immersed $2n$-submanifold
$F: M\ra N $ into a K\"{a}hler-Einstein manifold of complex dimension $2n$
and scalar curvature $R$.
We prove that, if  $M$ is compact, $n\geq 2$,  and $R<0$, then: 
$(i)$  Either F  has complex or Lagrangian directions;
$(ii)$ If $n=2$, M is oriented, and F has no complex directions, then 
it is a Lagrangian submanifold, generalizing the well-known case 
$n=1$ for minimal surfaces due to Wolfson.  We also prove that, 
 if $F$ has constant K\"{a}hler angles with no complex directions, and 
 is not  Lagrangian,  then $R=0$ must hold.
Our main tool is a formula on the Laplacian of a symmetric
function on the K\"{a}hler angles.\\[1mm] 
}
\markright{\sl\hfill Salavessa- Valli \hfill}
\section{Introduction}
\setcounter{Th}{0}
\setcounter{Pp}{0}
\setcounter{Cr} {0}
\setcounter{Lm} {0}
\setcounter{equation} {0}
Let $(N,J,g)$ be a K\"{a}hler manifold of complex dimension $2n$, with
complex structure $J$, Riemannian metric $g$, and
K\"{a}hler form $\omega(U,V)=g(JU,V)$. Let $F:M\ra N$ be an immersed
submanifold of real dimension $2n$. We take on $M$ its induced metric
$g_{M}=F^{*}g$.  The 2-form $F^{*}\omega$ on $M$, defines for
each $0\leq k\leq n$, the sets 
\[ \Omega_{2k}=\{p\in M: F^{*}\omega
\mbox{~has~ rank~} 2k \mbox{~at~the~point~}p\}.\] 
We will denote by $\Omega_{2k}^{0}$ the set of
interior points of $\Omega_{2k}$ in $M$. 
At each  point $p\in M$,  we identify
$F^{*}\omega$ with a skew-symmetric operator of $T_{p}M$ by using the 
musical isomorphism with respect to $g_{M}$, namely,\\[-6mm]
\[g_{M}(F^{*}\omega(X),Y)=F^{*}\omega(X,Y). \]
We take its polar 
decomposition
\begin{equation}
F^{*}\omega =  \tg \Jw
\end{equation}
where $J_{\omega}:T_{p}M\ra T_{p}M$ is a ( in fact unique) partial isometry with the same 
kernel ${\cal K}_{\omega}$ as  of $F^{*}w$, and where $\tg$ is the 
positive semidefinite 
operator $\tg = |F^{*}\omega|=\sqrt{-(F^{*}\omega)^{2} }$, with 
kernel ${\cal K}_{\omega}$ as well.
It turns out that
$J_{\omega}:{\cal K}^{\bot}_{\omega}\ra{\cal K}^{\bot}_{\omega}$ 
defines a complex structure  on ${\cal K}^{\bot}_{\omega}$, the orthogonal
complement of $\kw$ in $T_{p}M$. Moreover, $\Jw$ is 
$g_{M}$-orthogonal.  Since $\Fw$ 
is a normal operator, then $\Jw$ commutes with $\tg$. The map $P\ra |P|$ is Lipschitz in the space of normal operators of
a finite dimensional Hilbert space,  
for the Hilbert Schmidt norm ([Bh] p.215). Thus, the tensor
$\tg$ is continuous on all $M$ and locally Lipschitz. On 
 each $\Omega^{0}_{2k}$, $\kw$ and ${\cal K}^{\bot}_{\omega}$ are smooth sub-vector bundles of $TM$ and, from
the smoothness of the polar decomposition on invertible operators, 
  $\tg$ and $J_{\omega}$ 
are  smooth morphisms on these open sets. 
Let $\{ X_{\al},Y_{\al}\}_{1\leq \al \leq n}$ be a $g_{M}$-orthonormal  basis of $T_{p}M$, that diagonalizes 
$F^{*}\omega$ at $p$, that is 
\begin{equation}
F^{*}\omega=\bigoplus_{1\leq \al\leq n}  \left[ \begin{array}{cc}
             0 &  -a_{\al} \\
             a_{\al} &  0 \end{array} \right].
\end{equation}
Since $|a_{\al}|\leq 1$, then $a_{\al}=\cos\theta_{\al}$ for some 
angle $\theta_{\al}$.
We can reorder the diagonalizing basis such that 
 $\cos\theta_{1}\geq \cos\theta_{2}\geq \ldots \geq \cos\theta_{n}
\geq 0$. The angles  $\{\theta_{\al}\}_{ 1\leq \al\leq n}$
 are the \em K\"{a}hler angles  \em of $F$ at $p$. Thus, $\forall \al$, $F^{*}\omega(X_{\al})=\cos\theta_{\al}Y_{\al}$,
$F^{*}\omega(Y_{\al})=-\cos\theta_{\al}X_{\al}$ 
and if $k\geq 1$,
where $2k$ is the rank of $F^{*}\omega$ at $p$,
$J_{\omega}X_{\al}=Y_{\al}$ $~\forall \al \leq k$. The Weyl's
perturbation theorem (cf. [Bh]), 
 applied to the eigenvalues of the symmetric 
operator $|\Fw |$, shows that, ordering the $\cos\theta_{\al}$ in the
above way, the map 
$p\ra cos\theta_{\al}(p)$ is locally Lipschitz on $M$, for each $\al$.
Considering $\tg$
and $\Fw$  2-tensors, they are related by, 
\[\tg(X,Y)=\Fw(X,\Jw Y),~~~~ \forall X,Y\in T_{p}M.\]
Both $\tg$ and $\Fw$ are $\Jw$-invariant, that is 
$\tg(\Jw X,\Jw Y)=\tg (X,Y)$, $\Fw (\Jw X,\Jw Y)=\Fw (X,Y)$.
In particular, on $\Omega_{2n}$, $-\Fw$ is the K\"{a}hler form of $M$ with
respect the almost complex structure $\Jw$ and the Riemannian metric $\tg$.

A \em complex direction \em  of $F$ is a real two plane $P$ of $T_{p}M$
such that $dF(P)$ is a complex line of $T_{F(p)}N$, that is , 
$JdF(P)\subset dF(P)$. Similarly, $P$ is said to be a \em Lagrangian 
direction \em  of $F$ if $\omega$ vanishes on $dF(P)$,
that is, $JdF(P)\bot dF(P)$. The immersion $F$ has no complex directions  iff
$\cos\theta_{\al}<1$ $\forall \al$. The 2-plane $\{X_{\al},Y_{\al}\}$ is a complex direction of $M$ iff $J\circ dF = \pm dF \circ \Jw$ on that plane.
$M$ is a complex submanifold
iff  $cos\theta_{\al}=1$ $\forall \al$,  and is a 
Lagrangian submanifold iff  $cos\theta_{\al}=0$ $\forall \al$.
We say that $F$ has \em equal K\"{a}hler angles \em if $\theta_{\al}
=\theta$ $\forall \al$. Complex and Lagrangian submanifolds are 
examples of such case.

On each $\Omega_{2k}^{0}$, ${\cal K}_{\omega}^{\bot}$ 
is a smooth $J_{\omega}$-Hermitian
sub-vector bundle of $TM$.
Thus, for each $p_{0}\in \Omega_{2k}^{0}$ there exists a smooth local
$g_{M}$-orthonormal frame of ${\cal K}_{\omega}^{\bot}$ defined 
on a neigbourhood of  $p_{0}$, of the
form $X_{1},J_{\omega}X_{1},\ldots,X_{k},J_{\omega}X_{k}$. We may
enlarge it to a smooth $g_{M}$-orthonormal frame on $M$, on a 
neighbourhood of $p_{0}$
\begin{equation}
X_{1}\;, Y_{1}=J_{\omega}X_{1}\;,\ldots,X_{k}\;,Y_{k}=J_{\omega}X_{k}\;,
\,X_{k+1}\;,Y_{k+1}\;,\ldots\;,
X_{n}\;,Y_{n}\end{equation}
where $X_{k+1},Y_{k+1},\ldots X_{n},Y_{n}$ is any $g_{M}$-orthonormal
 frame of $\kw$. We may require it to be   a
diagonalizing basis of $F^{*}\omega$ at $p_{0}$. 
Note that in general it is not
possible to get smooth diagonalizing $g_{M}$-orthonormal frames in a whole
neighbourhood of a point $p_{0}$, unless $F^{*}\omega$ has distinct
non-zero eigenvalues $\cos\theta_{1},\ldots,\cos\theta_{n}$ at $p_{0}$,
or $F^{*}\omega$ has constant rank two, or $F$ has equal K\"{a}hler
angles.  

Let us denote by $\lnab{X}dF(Y)=\lnab{}dF(X,Y)$ the 
 second fundamental form of $F$.
It is a symmetric 2-tensor on $M$ that takes values on the normal bundle
$NM=(dF(TM))^{\bot}$. Let us denote by $(~)^{\bot}$ the orthogonal
projection of $F^{-1}TN$ onto the normal bundle.  We denote by $\lnab{}$
both Levi-Civita connections of $M$ and $N$.  We also denote by
$\lnab{}$ the induced connection on $F^{-1}TN$.
We take on $NM$ the usual connection $\lnabo{}$, given by $~\lnabo{X}U=
(\lnab{X}U)^{\bot}$, for $X\in T_{p}M$ and $U$ a smooth section of 
$NM\subset F^{-1}TN$.
$F$ is said to be \em totally geodesic \em if $\lnab{}dF =0$, and
$F$ is \em minimal \em if $trace_{g_{M}}\lnab{}dF =0$, that is $F$ is an harmonic immersion.

For a local frame  as in (1.3), and that diagonalizes $\Fw$ at $p_{0}$,
 we take the complex frame of $T^{c}M$
\[
Z_{\al}= \frac{X_{\al}-iY_{\al}}{2}= ``\al",~~~~~~~~ Z_{\bal}=
\frac{X_{\al}+iY_{\al}}{2}=``\bal",~~~~~~~~~\al\in\{1,\ldots,n\}.\] 
We extend by $\Co$-multilinearity  $g$, $g_{M}$,   $\tg$, $R^{N}$, $\Fw$, 
and any other tensors that may occur.  Sometimes we denote by $\langle,\rangle$ the $\Co{}$-bilinear extension of $g_{M}$.

If $p_{0}$ is a point {\bf without complex directions}, that is,
$\sin \theta_{\al}\neq 0~ \forall \al$, then,  $\{dF({\al}), dF({\bal}),$ $
U_{\al},U_{\bal}\}_{1\leq \al\leq n}$, where \\[-5mm]
\[
 U_{\al}=\frac{(JdF({\al}))^{\bot}}{\sin\theta_{\al}}~=~
\frac{JdF({\al})-i\cos\theta_{\al}dF({\al})}
{\sin\theta_{\al}},\]
 and $U_{\bal}=\overline{U_{\al}}$, is ( up to a constant
factor) an unitary basis of $T_{F(p_{0})}^{c}N$, with $U_{\al},U_{\bal}$
complex basis of the complexified normal bundle.
We denote by $R^{N}(U,V)W=-\lnab{U}\lnab{V}W+\lnab{V}\lnab{U}W+
\lnab{[U,V]}W$ and by $R^{N}(U,V,W,Z)=g(R^{N}(U,V)W,Z)$ the curvature 
tensor and Riemannian curvature tensor of $N$ respectively.
\begin{Lm}
At a point $p_{0}$ without complex directions, the Ricci tensor of $N$ is
given by, for $U,V\in T^{c}_{F(p_{0})}N$\\[-4mm] 
\[
Ricci^{N}(U,V)=\sum_{1\leq \mu\leq n}\frac{4}{\sin^{2}
\theta_{\mu}}R^{N}(U,JV,dF({\mu}),JdF({\bmu}) +i\cos\theta_{\mu}dF({\bmu}))
\]
\end{Lm}
\em Proof. \em Recall that
$R^{N}(JZ,JW,Z',W')=R^{N}(Z,W,JZ',JW')=R^{N}(Z,W,Z',W')$ and 
$Ricci^{N}(U,V)=-\frac{1}{2} trace \La{(} Z \ra R^{N}(U,JV)JZ \La{)}$.
Thus we have,\\[-5mm]
\begin{eqnarray*}
Ricci^{N}(U,V)&=&\sum_{\al} -\frac{1}{2}\La{(} 4 R^{N}(U,JV,JdF(\al),dF(\bal))
+4R^{N}(U,JV, JU_{\al},U_{\bal})\La{)}\\[-3mm]
\lefteqn{\!\!\!\!\!\!\!\!\!\!\!\!
\!\!\!\!\!\!\!\!\!\!\!\!\!\!\!\!\!\!\!\!\!\!\!\!\!\!\!\!\!\!\!\!
 =\!\!\sum_{\al}\! \frac{2}{\sin^{2}\theta_{\al}}
\LA{(}\!\! R^{N}(U,JV, dF(\al),\sin^{2}\theta_{\al} JdF(\bal))
\!+\!R^{N}(U,JV,dF(\al), JdF(\bal)\!+\!i\cos\theta_{\al}dF(\bal))}\\[-4mm]
&&-R^{N}(U,JV,dF(\al), -i\cos\theta_{\al}dF(\bal)-\cos^2\theta_{\al}JdF(\bal))
\LA{)} 
\end{eqnarray*}
leading the expression in the lemma. ~~~~\qed\\[4mm]
A first conclusion can be obtained from the above lemma:
\begin{Pp} Let $F:M\ra N$ be a totally geodesic immersion,  or more 
generally, let $F$ be an immersion such that the normal component of 
$F^{*}R^N(X,Y)Z$ vanishes, $\forall X, Y, Z \in TM$. 
If $F$ has no complex directions,
then $F^{*}\Psi=0$, where  $\Psi(U,V)=Ricci^{N}(JU,V)$,
$\forall U,V\in TN$, is the Ricci form of $N$. In particular, if $N$ is 
K\"{a}hler-Einstein of non-zero Ricci tensor, $F$ is a Lagrangian submanifold.
\end{Pp}
\em Proof. \em  We denote by $\lnab{X}\lnab{}dF$ the covariant derivative
of $\lnab{}dF$ as a section of $\bigotimes^{2}T^{*}M\otimes NM$.
Since $U_{\bmu}$ lies in the complexified normal bundle,
 using Codazzi's equation, $\forall X,Y\in
T_{p_{0}}M $
\[\begin{array}l
R^{N}(dF(X),dF(Y),dF({\mu}),JdF({\bmu}) +i\cos\theta_{\mu}dF({\bmu}))\\
=g\LA{(}-\lnab{X}\lnab{}dF (Y,\mu)+\lnab{Y}\lnab{}dF (X,\mu)\,,
\,JdF({\bmu}) +i\cos\theta_{\mu}dF({\bmu})\LA{)},\end{array}\]
that is zero for $F$ totally geodesic or
for $F$ with vanishing normal component of $ F^{*}R^N$.
 Therefore,  $ F^{*}\Psi(X,Y)=$
 $Ricci^{N}(JdF(X),dF(Y))$$=-Ricci^{N}(dF(X),JdF(Y))$ 
$=0$, by Lemma 1.1. ~~~~~\qed\\[4mm]
Henceforth, we assume $N$ is K\"{a}hler-Einstein with $Ricci^N = Rg$.
Our aim is to find conditions for a minimal immersion $F$ to be Lagrangian.
 Wolfson [W] proves that, for $n=1$,  if $F$ is a minimal real surface
without complex directions, immersed into a K\"{a}hler-Einstein surface of negative scalar curvature $R$, then $F$ is Lagrangian. His main tool is
a formula on the Laplacian of a smooth map $\kappa$ on the K\"{a}hler
angle of $F$, where the scalar curvature $R$ of $N$ appears. Here and in [S-V]
we generalize this map $\kappa$ to any dimension $n$ (see (2.2) below) and
compute its Laplacian. The main problem with this generalization for 
 $n\geq 2$, is that $\kappa$, although it is Lipschitz on $M$,
under the assumption of no existence of complex directions, it is no longer smooth on sets where $F$ might change the number of Lagrangian directions. Then, for instance,  we cannot use Stokes, or 
apply so easily the  maximum principle has Wolfson did.
The expression of $\triangle\kappa$ can be simplified in two cases.
The first  case is when  $F$ has equal K\"{a}hler angles, that we study
in [S-V]. In this case, we gain more regularity for $\kappa$, and we 
could conclude that, if $R\neq 0$ and $n=2$, then  either $F$ is a complex 
or a Lagrangian submanifold, and for $n\geq 3$, if $R<0$ and $F$ has no 
complex directions, then $F$ is Lagrangian (see [S-V]). 
The second case is when
$F$ is a  broadly-pluriminimal submanifold, 
a concept  we introduce in the next section. This is a concept that is
close to the product of minimal surfaces 
immersed into K\"{a}hler surfaces, possibly with different K\"{a}hler
angles. In this case the expression for $\triangle\kappa$
is the simpliest, and similar to the one of Wolfson [W]. 
Using this formula, we obtain a conclusion of Lagrangianity for the case $n=2$ and $R<0$ (Theorem 2.1), generalizing  the above result of Wolfson. 
\section{Broadly-pluriminimal submanifolds}
\setcounter{Th}{0}
\setcounter{Pp}{0}
\setcounter{Cr} {0}
\setcounter{Lm} {0}
\setcounter{equation} {0}
As in the previous section we let $(N,J,g)$  be a K\"{a}hler manifold of
complex dimension $2n$ and $F:M\ra N$ an immersed submanifold of real
dimension $2n$.\\[2mm]
{\bf Definition.} \em A map $F:M\ra N$ is said to be broadly-pluriminimal if\\[1mm]
$(i)$~~~~$F$ is minimal,\\
$(ii)$~~~On each  $\Omega_{2k}^{0}$, for  $ 1\leq k\leq n$,
 $F$ is pluriharmonic with respect to any local complex
structure $\tilde{J}=J_{\omega}\oplus J'$ where $J'$ is any
 $g_{M}$-orthogonal  complex structure
of ${\cal K}_{\omega}$. \em \\[4mm]
On the open set $\Omega_{2n}$, $(ii)$ means that $F$ is
pluriharmonic with respect to the complex structure $J_{\omega}$.
If $\kw=0$, 
we simply say that $F$ is \em pluriminimal. \em
We recall
that pluriharmonic maps are harmonic maps (see e.g [O-V]). 
If $F$ is broadly-pluriminimal, 
at each point $p\in \Omega_{2k}^{0}$, $k\geq 1$,  $\lnab{}dF$ is of type
$(2,0)+(0,2)$ for any complex structure $\tilde{J}$ of $T_{p}M$ of the
above form and properties,
 or equivalently, the $(1,1)$-part of it vanishes:
\\[-4mm]
\begin{equation}
 (\lnab{}dF)^{(1,1)}(X,Y)=\frac{1}{2}\La{(}\lnab{}dF(X,Y)
 +\lnab{}dF(\tilde{J}X,\tilde{J}Y)\La{)}=0.
\end{equation}\\
\em Examples. \em (1) Any minimal immersion of an oriented real surface
into a K\"{a}hler complex surface $F:M^{2}\ra N^{2}$ is broadly-pluriminimal.
In fact, we
have $F^{*}\omega= f Vol$, where $f:M\ra \R{}$ is a smooth map.  Any $g_{M}$-orthonormal basis $\{X,Y\}$,
diagonalizes $\Fw$, with $F^{*}\omega(X,Y)=\pm f=\pm \cos \theta$  
smooth everywhere.
In this case $J_{\omega}=\pm J_{M}$, wherever $f$ does not vanish, and where
$J_M$ is the na\-tu\-ral complex structure defined by a direct orthonormal basis. Since $F$ is harmonic and $n=1$, it is pluriharmonic with respect to $J_M$ and with respect to $-J_M$. 
 Condition $(ii)$ now follows, because, at each point
$p_{0}$ either ${\cal K}_{\omega}=T_{p_{0}}M$, that is $k=0$, or ${\cal
K}_{\omega}=0$, that is $k=n$.\\[2mm]
(2) The product, $F=F_{1}\times\ldots\times F_{n}:M=M_{1}\times\ldots \times M_{n}\ra N_{1}\times\ldots\times N_{n}$, 
of mi\-ni\-mal orientable surfaces
$F_{i}:M_{i}\ra N_{i}$ immersed into K\"{a}hler complex surfaces
and their reparametrizations, $\tilde{F}= F\circ
\phi:\tilde{M}\ra N_{1}\times\ldots\times N_{n}$,
where $\phi:\tilde{M}^{2n} \ra M= M_{1}\times\ldots \times M_{n}$ is a
diffeomorphism between manifolds,  are broadly-pluriminimal. More generaly,
the pro\-duct of broadly-pluriminimal submanifolds is a
 broadly-pluriminimal submanifold. We note that, to require the product
of two surfaces to have equal K\"{a}hler angles, implies each surface
to have constant K\"{a}hler angle (and equal to each other), for, the
two K\"{a}hler angles are indexed on independent variables. So, equal K\"{a}hler angles and  broadly-pluriminimality are independent concepts.
\\[2mm]
(3) If $F:M\ra N$ is a minimal Lagrangian submanifold, then $F$ is
broadly-pluriminimal. In this case ${\cal K}_{\omega}= TM$ and so
$\Omega^{0}_{2k}=\emptyset$ for $k\geq 1$. Thus $(ii)$ is
satisfied.\\[2mm]
(4) If $F:M \ra N$ is a complex submanifold, then $F$ is 
trivially pluriminimal. In this case $\Jw$ is the induced complex
structure from $J$.\\[2mm]
(5) If $F:M\ra N$ is a minimal  immersion with equal K\"{a}hler
angles and no complex directions, in [S-V] we prove that $F$ is
broadly-pluriminimal iff the isomorphism $\Phi : TM\ra NM$, given by
$~\Phi (X)= (JdF(X))^{\bot}$, is parallel, considering $TM$ with
the conformaly equivalent Riemannian metric 
$\hat{g}(X,Y)=g_{M}(X,Y)-g_{M}(\Fw (X), \Fw(Y) )=
\sin^{2}\theta g_{M}(X,Y)$. Moreover,
if this is the  case, for $n\geq 2$, it turns out that the K\"{a}hler angle is constant, and so, $\Phi: (TM,g_{M}) \ra (NM,g)$ is a parallel homothetic isomorphism. Furthermore, if $F$ is not Lagrangian, $N$ must be Ricci-flat, 
as a consequence 
of Corollary 2.1 given in the next section of this paper. This example shows how natural is the definition of broadly-pluriminimal.\\[2mm]
(6) Any minimal immersion $F:M^{2n}\ra T^{2n}$ into the flat complex
torus with no Lagrangian directions and such that  $(M,\Jw,g_{M})$
is K\"{a}hler, is a pluriminimal submanifold. In fact, from Gauss
equation, $\sum_{\al,\mu}R^{M}(\mu,\al,\bmu,\bal)=-\sum_{\al,\mu}
\|\lnab{}dF(\al,\bmu)\|^{2}$, where $R^{M}$ is the Riemannian curvature 
tensor of $M$, that is of type $(1,1)$ with respect to  $\Jw$.\\[2mm]
(7) Let $(N,I,J,g)$ be a hyper-K\"{a}hler manifold of real dimension 8, 
and $F:M\ra N$  a minimal immersion of a 4-dimensional 
sub\-ma\-ni\-fold  with non-negative isotropic sectional curvature, and  such that $\forall \nu\phi\in S^{2}$, $F$ has equal K\"{a}hler angles with respect 
to the complex structure $J_{\nu\phi}=\cos\nu I+\sin\nu\cos\phi J
+\sin\nu\sin\phi K$, where $K=IJ$. Set for each unit vector $X\in TM$,
$H_{X}=span\{ X, IX, JX, KX\}$.
In [S-V] we prove that,
if $\exists p\in M$ and  $\exists X\in T_{p}M$, unit vector,
such that $dim(T_{p}M\cap H_X)\geq 2$, then there exists 
$\nu\phi\in S^{2}$ such that $M$ is a $J_{\nu\phi}$-complex 
submanifold. Furthermore, if $J_{\nu\phi}=I$ then $F: M \ra (N,I,g)$ is obviously pluriminimal. If $J_{\nu\phi}\neq I$ but
 $T_{p}M\cap H_X^{\bot}\neq \{ 0\}$, 
then  $F: M \ra (N,I,g)$ is still pluriminimal, with constant
K\"{a}hler angle $\nu$ or $\nu + \pi$.\\[2mm]

If $\{\theta_{\al}\}_{1\leq n\leq n}$ 
are the  K\"{a}hler angles of $F$, 
$~g_{M}\pm\tg$, where $\tg$ is given in (1.1), is represented 
in the unitary basis $\{\sqrt{2}\al,\sqrt{2}\bal\}$ of $T_{p}M$, for $p$
near $p_{0}$, 
by a $2n\times 2n $ matrix that at $p_{0}$ is the diagonal  matrix $D(1\pm\cos\theta_{1},\ldots,1\pm\cos\theta_{n},
1\pm\cos\theta_{1},\ldots,1\pm\cos\theta_{n})$. Thus, 
$ det (g_{M}\pm\tg)=\prod_{1\leq\al\leq n}(1\pm\cos\theta_{\al})^{2}$.
If $p_{0}$ is a point without complex directions, $\cos\theta_{\al}\neq 1$,
$\forall\al\in\{1,\ldots,n\}$, and so $\tg<g_{M}$. Thus, is $F$
 has no complex directions we may consider the map
\begin{equation}
\kappa =\frac{1}{2}\log \left( \frac{det(g_{M}+\tg)}{det(g_{M}-\tg)}\right)
=\sum_{1\leq \al\leq
n}\log\LA{(}\frac{1+\cos\theta_{\al}}{1-\cos\theta_{\al}}\LA{)}. 
\end{equation}
This map $\kappa$ is non-negative, continuous on $M$ and smooth on each
$\Omega_{2k}^{0}$. It is an increasing map on each $\cos\theta_{\al}$. 
In [S-V] we compute $\triangle \kappa$ on $\Omega_{2k}^{0}$ ( see Proposition
2.1 below), for any minimal immersion $F$. In this paper we  will compute $\triangle \kappa$ for the case of $F$ broadly-pluriminimal (Proposition
2.2). Then, using
such formula we can prove the statements given in the abstract, 
and some few other ones.
\subsection{The computation of $\triangle\kappa$}
Let $\{X_{\al}, Y_{\al}\}$ be a local $g_{M}$-orthonormal frame 
satisfying the conditions in (1.3), and that  diagonalizes 
$\Fw$ at $p_{0}$.  We define a local $g_{M}$-orthogonal complex structure
on a neighbourhood of $p_{0}\in \Omega^{0}_{2k}$ as
$\tilde{J}=J_{\omega}\oplus J'$, where $J_{\omega}$ is defined on
${\cal K}^{\bot}_{\omega}$ and $J'$ is the local complex structure on
${\cal K}_{\omega}$, defined on a neigbourhood of $p_{0}$ by
 $J'Z_{\al}=iZ_{\al}$, $J'Z_{\bal}=-iZ_{\bal}$, $\forall \al\geq k+1$.
On a neigbourhood of $p_{0}$, $Z_{\al}$ is of type (1,0) with respect to
$\tilde{J}$,  $\forall\al $, and  
$Z_{\al}$ and $Z_{\bal}$ are in ${\cal K}_{\omega}^{c}$ 
~$\forall \al \geq k+1$. Since $\tilde{J}$ is $g_{M}$-orthogonal,
 $\forall \al,\be$, on a neigbourhood of $p_{0}$,\\[-4mm] 
\begin{equation}
\langle \lnab{Z}\tilde{J}({\al}),{\be}\rangle =2i\langle\lnab{Z} {\al},
{\be}\rangle= -\langle
{\al},\lnab{Z}\tilde{J}({\be})\rangle,~~~~~~~~~~ \langle
\lnab{Z}\tilde{J}({\al}),{\bbe}\rangle =0. 
\end{equation}
In particular $\lnab{Z}\tilde{J}({\al})$ is of type $(0,1)$.
Note that, considering $\tg$ and $\Fw$  2-tensors, $\tilde{g}(X,Y)=\Fw(X,\tilde{J}Y)$
still holds $\forall X,Y\in T_{p}M$.
Set $\tg_{AB}=\tg({A},{B})$, and define 
$\overline{\bar{B}}=B$, $\forall A,B\in \{1,\ldots, n, \bar{1},
\ldots, \bar{n}\}$.  Let $\epsilon_{\al}=1$,
 $\epsilon_{\bal}=-1$ $\forall \al\in\{1,\ldots,n\}$.
Then $\forall 1\leq \al,\be\leq n$,  $\forall A,B \in \{1,\ldots, n, \bar{1},
\ldots, \bar{n}\}$ and $\forall C\in\{1,\ldots,n\}\cup \{\overline{k+1},\ldots,\bar{n}\}$,
\begin{equation}\left.\begin{array}{ll}
F^{*}\omega({\al},{C})=g(JdF({\al}),dF(C))=0 &~~~~~~\forall p
\mbox{~near~} p_{0}\\[1mm]
F^{*}\omega({\al},{\bbe})=g(JdF({\al}),dF({\bbe}))
=\frac{i}{2}\delta_{\al\be}\cos\theta_{\al} &~~~~~~\mbox{at~}p_{0}\\[1mm]
\tg_{AB}=i\epsilon_{B}F^{*}\omega({A},{B})=i\epsilon_{B}g(JdF({A}),dF({B}))
&~~~~~~\forall p \mbox{~near~} p_{0}\\[1mm] 
\tg_{\al C}=\tg_{\bal\bar{C}}=0
&~~~~~~\forall p \mbox{~near~} p_{0}\\[1mm]
\tg_{\al\bbe}=\tg_{\bal\be}=\frac{1}{2}\delta_{\al\be}\cos\theta_{\al}
&~~~~~~\mb
{at~}p_{0} \end{array}\right\}\end{equation}
Now we compute the covariant derivative of $F^{*}\omega$. Let $p\in M$, $X,Y,Z\in T_{p}M$. Then
\begin{equation}\begin{array}{ccc}
d(g(JdF(X),dF(Y)))(Z)\!\!&\!=\!&
g(J\lnab{Z}dF(X),dF(Y))+g(JdF(\lnab{Z}X),dF(Y))\nonumber\\[1mm]
\!\!\!\!\!\!\!&&
+g(JdF(X),\lnab{Z}dF(Y))+g(JdF(X),dF(\lnab{Z}Y))\end{array}
 \end{equation}\\[-6mm]
And so\\[-6mm]
\begin{equation} \lnab{Z}F^{*}\omega
(X,Y)=-g(\lnab{Z}dF(X),JdF(Y))+g(\lnab{Z}dF(Y),JdF(X)). \end{equation}
For simplicity of notation we denote by\\[-6mm]
\[ \gf{X}{Y}{Z}=\gdf{X}{Y}{Z}. \]
 From (2.5) and (2.4) we have
\begin{Lm}
$\forall p$ near $p_{0}\in\Omega^{0}_{2k}$, $Z\in T^{c}_{p}M$, and
$\mu,\ga\in\{1,\ldots,n\}$
\begin{eqnarray*}
d\tg_{\mu\bga}(Z) & = & i\gf{Z}{\mu}{\bga}-i\gf{Z}{\bga}{\mu}
+2\sum_{ \rho}\la{(}\langle
\lnab{Z}{\mu},{\brho}\rangle\tg_{\rho\bga} +\langle
\lnab{Z}{\bga},{\rho}\rangle\tg_{\mu\brho}\la{)}\\[-2mm]
0=d\tg_{\mu\ga}(Z)	& =
&-i\gf{Z}{\mu}{\ga}+i\gf{Z}{\ga}{\mu} +2\sum_{ \rho}\la{(}
\langle\lnab{Z}{\mu},{\rho}\rangle\tg_{\brho\ga} -\langle
\lnab{Z}{\ga},{\rho}\rangle\tg_{\mu\brho}\la{)}. \end{eqnarray*}\\[-7mm]
In particular, at $p_{0}$\\[-9mm]
\begin{eqnarray*}
d\tg_{\mu\bga}(Z) & = &i\gf{Z}{\mu}{\bga}-i\gf{Z}{\bga}{\mu}
 -(\cos\theta_{\mu} -\cos\theta_{\ga})\langle
\lnab{Z}{\mu},{\bga}\rangle\\[-1mm]
0=d\tg_{\mu\ga}(Z)	& = & -i\gf{Z}{\mu}{\ga}+i\gf{Z}{\ga}{\mu}
+(\cos\theta_{\mu} +\cos\theta_{\ga})\langle \lnab{Z}{\mu},{\ga}\rangle.
\end{eqnarray*}
\end{Lm}
\begin{Pp} {\bf [S-V]}
~If F is minimal without complex directions, and  $\{X_{\al},Y_{\al}\}$ is a local orthonormal  frame of the form $(1.3)$ that diagonalizes $\Fw$ at
 $p_{0}\in\Omega^{0}_{2k}$, $0\leq k\leq n$, then at $p_{0}$
\begin{eqnarray}
\triangle \kappa
&=&4i\sum_{\be} Ricci^{N}(JdF(\be),dF(\bbe))\non\\[-3mm]
&&+\sum_{\be,\mu}\!\frac{32}{\sin^{2}\theta_{\mu}} 
Im \La{(}R^{N}(dF(\be),dF(\mu),dF(\bbe),
JdF(\bmu)\!+\!i\cos\theta_{\mu}dF(\bmu))\La{)}\non\\[-2mm]
&&-\sum_{\be,\mu,\rho}\frac{64(\cos\theta_{\mu}
\!+\!\cos\theta_{\rho})}{\sin^{2}\theta_{\mu}\sin^{2}\theta_{\rho}}
Re\la{(}\gf{\be}{\mu}{\brho}\gf{\bbe}{\rho}{\bmu}\la{)}\non\\
&& +\sum_{\be,\mu,\rho}
\frac{32(\cos\theta_{\rho}-\cos\theta_{\mu})}
{\sin^{2}\theta_{\mu}\sin^{2}\theta_{\rho}}\;(|\gf{\be}{\mu}{\rho}|^{2}
+|\gf{\bbe}{\mu}{\rho}|^{2})\\
&&+\sum_{\be,\mu,\rho}\frac{32(\cos\theta_{\mu}+\cos\theta_{\rho})}
{\sin^{2}\theta_{\mu}}\,\LA{(}
|\langle\lnab{\be}\mu,\rho\rangle|^{2} +
|\langle\lnab{\bbe}\mu,\rho\rangle|^{2}\LA{)}.\non
\end{eqnarray}
\end{Pp}
Now we get
\begin{Pp} If $F$ is broadly-pluriminimal without complex directions,
and  $\{X_{\al},Y_{\al}\}$ is a local orthonormal 
 frame of the form $(1.3)$ that diagonalizes $\Fw$ at
 $p_{0}\in\Omega^{0}_{2k}$, $0\leq k\leq n$, then at $p_{0}$
\begin{equation}
\triangle
\kappa = 4i \sum_{1\leq\be\leq n} Ricci^{N}(JdF({\be}),dF({\bbe}))
\end{equation}
\end{Pp}
\em Proof. \em ~ First, we rewrite $\triangle\kappa$ of Proposition 2.1
 in terms of
$(\lnab{}dF)^{(1,1)}$, the symmetric tensor given by
(2.1)  on  a neighbourhood of $p_{0}$, and it is 
 the (1,1)-part of $\lnab{}dF$ with respect to the complex structure
$\tilde{J}$. From Lemma 2.1, at $p_{0}$,
\[|\gf{Z}{\mu}{\rho}|^{2}=|\gf{Z}{\rho}{\mu}|^{2} +(\cos\theta_{\mu}
+\cos\theta_{\rho})^{2}|\langle\lnab{Z}\mu, \rho\rangle|^{2}
+2(\cos\theta_{\mu}+\cos\theta_{\rho}) Im( \langle \lnab{Z}\mu,
\rho\rangle\gf{\bar{Z}}{\brho}{\bmu}).\]
Hence, in the expression of $\triangle \kappa$ of Proposition 2.1,
\begin{eqnarray}
(2.7)&=&\sum_{\be,\mu,\rho}
\frac{16(\cos\theta_{\rho}\!-\!\cos\theta_{\mu})}
{\sin^{2}\theta_{\mu}\sin^{2}\theta_{\rho}}\la{(}
|\gf{\be}{\mu}{\rho}|^{2}-|\gf{\be}{\rho}{\mu}|^{2}+
|\gf{\bbe}{\mu}{\rho}|^{2}-|\gf{\bbe}{\rho}{\mu}|^{2}\la{)}\nonumber\\
&=&\sum_{\be,\mu,\rho}
\frac{16(\cos\theta_{\rho}\!-\!\cos\theta_{\mu})
(\cos\theta_{\mu}\!+\!\cos\theta_{\rho})^{2}}
{\sin^{2}\theta_{\mu}\sin^{2}\theta_{\rho}}\la{(}
|\langle\lnab{\be}\mu,\rho\rangle|^{2}
+|\langle\lnab{\bbe}\mu,\rho\rangle|^{2}\la{)}\nonumber\\[-2mm]
&&+\sum_{\be,\mu,\rho}
\frac{32(\cos\theta_{\rho}\!-\!\cos\theta_{\mu})
(\cos\theta_{\mu}\!+\!\cos\theta_{\rho})}
{\sin^{2}\theta_{\mu}\sin^{2}\theta_{\rho}}Im \la{(}\langle\lnab{\be}
\mu,\rho\rangle\gf{\bbe}{\brho}{\bmu}+\langle\lnab{\bbe}
\mu,\rho\rangle\gf{\be}{\brho}{\bmu}\la{)}\nonumber\\
&=&\sum_{\be,\mu,\rho}
\frac{16(\sin^{2}\theta_{\mu}\!-\!\sin^{2}\theta_{\rho})
(\cos\theta_{\mu}\!+\!\cos\theta_{\rho})}
{\sin^{2}\theta_{\mu}\sin^{2}\theta_{\rho}}\la{(}
|\langle\lnab{\be}\mu,\rho\rangle|^{2}
+|\langle\lnab{\bbe}\mu,\rho\rangle|^{2}\la{)}\nonumber\\[-2mm]
&&+\sum_{\be,\mu,\rho}\!
\frac{32(\sin^{2}\theta_{\mu}\!-\!\sin^{2}\theta_{\rho})}
{\sin^{2}\theta_{\mu}\sin^{2}\theta_{\rho}}Im \la{(}\langle\lnab{\be}
\mu,\rho\rangle\gf{\bbe}{\brho}{\bmu}+\langle\lnab{\bbe}
\mu,\rho\rangle\gf{\be}{\brho}{\bmu}\la{)}\nonumber\\
&=&\!\!\!\!\sum_{\be,\mu,\rho}\!16
\La{(}\frac{1}{\sin^{2}\theta_{\rho}}-\frac{1}{\sin^{2}\theta_{\mu}}\La{)}
(\cos\theta_{\mu}\!+\!\cos\theta_{\rho})
\la{(}|\langle\lnab{\be}\mu,\rho\rangle|^{2}
\!+\!|\langle\lnab{\bbe}\mu,\rho\rangle|^{2}\la{)}\\[-2mm]
&&+\sum_{\be,\mu,\rho}32
\La{(}\frac{1}{\sin^{2}\theta_{\rho}}-\frac{1}{\sin^{2}\theta_{\mu}}\La{)}
Im\LA{(}\langle\lnab{\be}\mu,\rho
\rangle\gf{\bbe}{\brho}{\bmu}+\langle\lnab{\bbe}\mu,\rho
\rangle\gf{\be}{\brho}{\bmu}\LA{)}.\nonumber
\end{eqnarray}
Note that $(2.9)=0$, because it is the product of skew-symmetric factor 
on $\rho,\mu$ with a symmetric one. Hence, interchanging $\mu$ with $\rho$
when necessary,
\begin{eqnarray*}
(2.7) &=&\sum_{\be,\mu,\rho}32
\La{(}\frac{1}{\sin^{2}\theta_{\rho}}-\frac{1}{\sin^{2}\theta_{\mu}}\La{)}
Im \la{(}\langle\lnab{\be}\mu,\rho
\rangle\gf{\bbe}{\brho}{\bmu}+\langle\lnab{\bbe}\mu,\rho
\rangle\gf{\be}{\brho}{\bmu}\la{)}\\[-2mm]
&=& \sum_{\be,\mu,\rho}
\frac{32}{\sin^{2}\theta_{\mu}}
Im \La{(}\langle\lnab{\be}\rho,\mu
\rangle\gf{\bbe}{\bmu}{\brho}+\langle\lnab{\bbe}\rho,\mu
\rangle\gf{\be}{\bmu}{\brho}-
\langle\lnab{\be}\mu,\rho
\rangle\gf{\bbe}{\brho}{\bmu}-\langle\lnab{\bbe}\mu,\rho
\rangle\gf{\be}{\brho}{\bmu}\La{)}.
\end{eqnarray*}
Therefore, and since $\langle \lnab{Z}\rho , \mu\rangle =-
\langle \lnab{Z}\mu , \rho \rangle$,
\begin{eqnarray}
\triangle \kappa &=&
4i\sum_{\be} Ricci^{N}(JdF(\be),dF(\bbe))\nonumber\\[-3mm]
 &&+\sum_{\be,\mu}\!\frac{32}{\sin^{2}\theta_{\mu}} 
Im \La{(}R^{N}(dF(\be),dF(\mu),dF(\bbe),
JdF(\bmu)+i\cos\theta_{\mu}dF(\bmu))\La{)}\\[-1mm]
&&-\sum_{\be,\mu,\rho}\frac{64(\cos\theta_{\mu}
\!+\!\cos\theta_{\rho})}{\sin^{2}\theta_{\mu}\sin^{2}\theta_{\rho}}
Re\la{(}\gf{\be}{\mu}{\brho}\gf{\bbe}{\rho}{\bmu}\la{)}\nonumber\\[-1mm]
&&-\sum_{\be,\mu,\rho}
\frac{32}{\sin^{2}\theta_{\mu}}
Im\la{(}\langle\lnab{\be}\mu,\rho\rangle\gf{\bbe}{\bmu}{\brho}
+\langle\lnab{\be}\mu,\rho\rangle\gf{\bbe}{\brho}{\bmu}
 +\langle\lnab{\bbe}\mu,\rho\rangle\gf{\be}{\bmu}{\brho}
+\langle\lnab{\bbe}\mu,\rho\rangle\gf{\be}{\brho}{\bmu}
\la{)}\nonumber\\[-1mm]
&&+\sum_{\be,\mu,\rho}\frac{32(\cos\theta_{\mu}+\cos\theta_{\rho})}
{\sin^{2}\theta_{\mu}}\,\LA{(}
|\langle\lnab{\be}\mu,\rho\rangle|^{2} +
|\langle\lnab{\bbe}\mu,\rho\rangle|^{2}\LA{)}.\nonumber
\end{eqnarray}
 Derivating (2.1), considering $\lnab{}dF$ and $
(\lnab{}dF)^{(1,1)}$ both with values in the normal bundle, we get\\[-9mm]
\begin{eqnarray*}
\mbox{\Large $\nabla$}_{Z}(\lnab{}dF)^{(1,1)}(X,Y)&=&\frac{1}{2}
\LA{(}\mbox{\Large $\nabla$}_{Z}\lnab{}dF(X,Y)+\mbox{\Large
$\nabla$}_{Z}\lnab{}dF(\tilde{J}X,\tilde{J}Y)\\[-3mm]
&&~~~~+\lnab{}dF(\lnab{Z}\tilde{J}(X),\tilde{J}Y)+
\lnab{}dF(\tilde{J}X,\lnab{Z}\tilde{J}(Y))\LA{)}.
\end{eqnarray*}\\[-7mm]
Since $\tilde{J}Z_{\mu}=iZ_{\mu}$ and $\tilde{J}Z_{\bmu}=-iZ_{\bmu}$
$\forall \mu$, on a neighbourhood of $p_{0}$, we have, $\forall \al,\be$
\begin{eqnarray}
\mbox{\Large $\nabla$}_{Z}(\lnab{}dF)^{(1,1)}({\al},{\be})
&=&\frac{i}{2}\La{(}\lnab{}dF(\lnab{Z}\tilde{J}({\al}),{\be})+
\lnab{}dF({\al},\lnab{Z}\tilde{J}({\be}))\La{)},\\ [2mm]
\mbox{\Large
$\nabla$}_{Z}(\lnab{}dF)^{(1,1)}({\al},{\bbe})
&=& \mbox{\Large $\nabla$}_{Z}\lnab{}dF({\al},{\bbe})\non\\[-2mm]
&&+\frac{i}{2}\La{(}-\lnab{}dF(\lnab{Z}\tilde{J}({\al}),{\bbe})+
\lnab{}dF({\al},\lnab{Z}\tilde{J}({\bbe}))\La{)}.
\end{eqnarray}
Since $F$ is minimal,\\[-5mm]
\[\sum_{\be}\mbox{\Large $\nabla$}_{Z}\lnab{}dF({\be},{\bbe})
=\frac{1}{4}trace_{g_{M}} \mbox{\Large $\nabla$}_{Z}\lnab{}dF
=\frac{1}{4}\mbox{\Large $\nabla$}_{Z}(trace_{g_{M}}\lnab{}dF)=0.\]\\[-5mm]
Then applying Codazzi's equation and noting that
$JdF({\bmu})+i\cos\theta_{\mu}dF({\bmu})$ is in the complexified normal
bundle,
\begin{equation}
\!\!\sum_{\be}\!R^{N}(\!dF(\!{\be}),\!dF(\!{\mu}\!),dF({\bbe}),
\!Jd\!F({\bmu})\!\!+\!i\!\cos\theta_{\mu}dF(\bmu))
=\!\sum_{\be}\! g\la{(}\!\!-\!\mbox{\Large $\nabla$}_{\!{\be}}\!
\lnab{}dF({\mu},{\bbe}),JdF({\bmu})\la{)}.
\end{equation}
Using (2.3), (2.12) and (2.13)
\begin{eqnarray*}
(2.10)
&=&-\sum_{\be,\mu}\frac{32}{\sin^{2}\theta_{\mu}}Im\La{(} g\la{(}\mbox{\Large$\nabla$}_{\be}(\lnab{}dF)^{(1,1)}(\mu,\bbe),
JdF(\bmu))\La{)}\\[-2mm]
&& +\sum_{\be,\mu,\rho}\frac{64}{\sin^{2}\theta_{\mu}}Im\La{(} \langle\lnab{\be}\mu,\rho\rangle
\gf{\brho}{\bbe}{\bmu}\La{)} +\sum_{\be,\mu,\rho}\frac{64}{\sin^{2}\theta_{\mu}}Im\La{(}
\langle\lnab{\be}\bbe,\brho\rangle\gf{\mu}{\rho}{\bmu}\La{)}.
\end{eqnarray*}
Consequently,
\begin{eqnarray}
\triangle \kappa
&=&4i\sum_{\be} Ricci^{N}(JdF(\be),dF(\bbe))
-\sum_{\be,\mu}\frac{32}{\sin^{2}\theta_{\mu}}Im\La{(} g\la{(}\mbox{\Large$\nabla$}_{\be}(\lnab{}dF)^{(1,1)}(\mu,\bbe),
JdF(\bmu))\La{)}\nonumber\\[-2mm]
&& +\sum_{\be,\mu,\rho}\frac{64}{\sin^{2}\theta_{\mu}}Im\La{(} \langle\lnab{\be}\mu,\rho\rangle
\gf{\brho}{\bbe}{\bmu}\La{)}\\[-2mm]
&& +\sum_{\be,\mu,\rho}\frac{64}{\sin^{2}\theta_{\mu}}Im\La{(}
\langle\lnab{\be}\bbe,\brho\rangle\gf{\mu}{\rho}{\bmu}\La{)}
-\sum_{\be,\mu,\rho}\!\!\!\frac{64(\cos\theta_{\mu}
\!+\!\cos\theta_{\rho})}{\sin^{2}\theta_{\mu}\sin^{2}\theta_{\rho}}
Re\la{(}\gf{\be}{\mu}{\brho}\gf{\bbe}{\rho}{\bmu}\!\la{)}\nonumber\\[-2mm]
&&-\sum_{\be,\mu,\rho}\frac{32}{\sin^{2}\theta_{\mu}}
Im \la{(}\langle\lnab{\be}\mu,\rho
\rangle\gf{\bbe}{\bmu}{\brho}
+\langle\lnab{\be}\mu,\rho
\rangle\gf{\bbe}{\brho}{\bmu}\la{)}\\[-2mm]
&&-\sum_{\be,\mu,\rho}\frac{32}{\sin^{2}\theta_{\mu}}
Im \la{(} \langle\lnab{\bbe}\mu,\rho
\rangle\gf{\be}{\bmu}{\brho}
+\langle\lnab{\bbe}\mu,\rho
\rangle\gf{\be}{\brho}{\bmu}\la{)}\\[-2mm]
&&+\sum_{\be,\mu,\rho}\frac{32(\cos\theta_{\mu}+\cos\theta_{\rho})}
{\sin^{2}\theta_{\mu}}\,\LA{(}
|\langle\lnab{\be}\mu,\rho\rangle|^{2} +
|\langle\lnab{\bbe}\mu,\rho\rangle|^{2}\LA{)}.
\end{eqnarray}
By lemma 2.1\\[-5mm]
\[
(2.14)+(2.15)=\!\!
\sum_{\be,\mu,\rho}\!\frac{32}{\sin^{2}\theta_{\mu}}Im\La{(}
\langle\lnab{\be}\mu,\rho\rangle (\gf{\bbe}{\brho}{\bmu}-
\gf{\bbe}{\bmu}{\brho})\La{)}
=-\!\!\sum_{\be,\mu,\rho}\!\!\frac{32(\cos\theta_{\mu}\!+\!\cos\theta_{\rho})}
{\sin^{2}\theta_{\mu}}|\langle\lnab{\be}\mu,\rho\rangle |^{2}
\]
that cancels with some term of (2.17). Similarly\\[-6mm]
\begin{eqnarray*}
\lefteqn{(2.16)+
\sum_{\be,\mu,\rho}\frac{32(\cos\theta_{\mu}+\cos\theta_{\rho})}
{\sin^{2}\theta_{\mu}}|\langle\lnab{\bbe}\mu,\rho\rangle |^{2}=}\\[-2mm]
&=& -\sum_{\be,\mu,\rho}\frac{32}{\sin^{2}\theta_{\mu}}Im\LA{(}
\langle\lnab{\bbe}\mu,\rho\rangle \la{(}
\gf{\be}{\bmu}{\brho}-\gf{\be}{\brho}{\bmu}\la{)}\LA{)}
-\sum_{\be,\mu,\rho}\frac{64}{\sin^{2}\theta_{\mu}}Im\La{(}
\langle\lnab{\bbe}\mu,\rho\rangle \gf{\be}{\brho}{\bmu}\La{)}\\
&&+\sum_{\be,\mu,\rho}\frac{32(\cos\theta_{\mu}+\cos\theta_{\rho})}
{\sin^{2}\theta_{\mu}}|\langle\lnab{\bbe}\mu,\rho\rangle |^{2}\\
&=& -\sum_{\be,\mu,\rho}\frac{64}{\sin^{2}\theta_{\mu}}Im\La{(}
\langle\lnab{\bbe}\mu,\rho\rangle \gf{\be}{\brho}{\bmu}\La{)}
\end{eqnarray*}
Then, 
\begin{eqnarray}
\triangle \kappa
&=&4i\sum_{\be} Ricci^{N}(JdF(\be),dF(\bbe))
-\sum_{\be,\mu}\frac{32}{\sin^{2}\theta_{\mu}}Im\La{(} g\la{(}\mbox{\Large$\nabla$}_{\be}(\lnab{}dF)^{(1,1)}(\mu,\bbe),
JdF(\bmu)\la{)}\La{)}\nonumber\\[-2mm]
&&-\sum_{\be,\mu,\rho}\!\!\!\frac{64(\cos\theta_{\mu}
\!+\!\cos\theta_{\rho})}{\sin^{2}\theta_{\mu}\sin^{2}\theta_{\mu}}
Re\la{(}\gf{\be}{\mu}{\brho}\gf{\bbe}{\rho}{\bmu}\!\la{)}\nonumber\\[-2mm]
&& +\sum_{\be,\mu,\rho}\frac{64}{\sin^{2}\theta_{\mu}}Im\La{(}
\langle\lnab{\be}\bbe,\brho\rangle\gf{\rho}{\mu}{\bmu}\La{)}\\[-2mm]
&&-\sum_{\be,\mu,\rho}\frac{64}{\sin^{2}\theta_{\mu}}
Im \la{(}\langle\lnab{\bbe}\mu,\rho
\rangle\gf{\be}{\brho}{\bmu}\la{)}
\end{eqnarray}
Now, from (2.3), $(\lnab{\be}\bbe)^{(1,0)}=\frac{i}{2}\lnab{\be}\tilde{J}(\bbe)$, and so
\[
(2.18)=\sum_{\be,\mu,\rho}\frac{32}{\sin^{2}\theta_{\mu}}Im\La{(}
g(\lnab{}dF( \frac{i}{2}\lnab{\be}\tilde{J}(\bbe),\mu),JdF(\bmu)\La{)}\]
and from (2.11), (2.3), and that $\lnab{\bbe}\tilde{J}(\mu)$ is of type
$(0,1)$
\begin{eqnarray*}
(2.19)&=&-\sum_{\be,\mu,\rho}\frac{32}{\sin^{2}\theta_{\mu}}Im\La{(}
g\la{(}\lnab{}dF(\be,-\frac{i}{2}\lnab{\bbe}
\tilde{J}(\mu)),JdF(\bmu)\la{)}\La{)}\\[-2mm]
&=&\sum_{\be,\mu,\rho}\frac{32}{\sin^{2}\theta_{\mu}}Im\La{(}
g\la{(}\lnab{\bbe}(\lnab{}dF)^{(1,1)}( \be,\mu),JdF(\bmu)\la{)}\La{)}\\[-2mm]
&&-\!\!\!\!\sum_{\be,\mu,\rho}\frac{32}{\sin^{2}\theta_{\mu}}Im\La{(}
g\la{(}\lnab{}dF(\frac{i}{2}\lnab{\bbe}
\tilde{J}(\be),\mu),JdF(\bmu)\la{)}\La{)}
\end{eqnarray*}
Using the unitary basis $\{\sqrt{2}\al,\sqrt{2}\bal\}$ of $T_{p}M$, for $p$
near $p_{0}$, $g_{M}+\tg$ is represented by the matrix 
\[g_{M}\pm\tg=\left[ \begin{array}{cc}
\delta_{\al\ga}\pm2\tg_{\al\bga} & 0\\
0& \delta_{\al\ga}\pm2\tg_{\bal\ga},\end{array}\right]\]
with $\tg_{\mu\brho}=\tg_{\bmu\rho}$.
 This matrix is at the point $p_{0}$ the diagonal matrix
$D(1\pm \cos\theta_{1}, \ldots , 1\pm\cos\theta_{n}, 1\pm \cos\theta_{1}, \ldots , 1\pm\cos\theta_{n})$.  Thus (cf. lemma 5.2 of [S-V]), $\forall
Z\in T_{p_{0}}M$,
\[d(det(g_{M}\pm \tg ))(Z)=\pm 4\sum_{\mu}\frac{det(g_{M}\pm\tg)}{(1\pm\cos\theta_{\mu})}d\tg_{\mu\bmu}(Z).\]
Then, using Lemma 2.1\\[-6mm]
\begin{eqnarray}
2d\kappa_{p_{0}}(Z)&=& \frac{d(det(g_{M}+\tg))(Z)}{det(g_{M}+\tg)}
-\frac{d(det(g_{M}-\tg))(Z)}{det(g_{M}-\tg)}\non\\[-1mm]
&=& 4\sum_{\mu}\frac{1}{(1+\cos\theta_{\mu})}
d\tg_{\mu\bmu}(Z)+4\sum_{\mu}\frac{1}{(1-\cos\theta_{\mu})}
d\tg_{\mu\bmu}(Z)\\[-1mm]  
&=&8\sum_{\mu}\frac{i}{\sin^{2}\theta_{\mu}}\La{(}\langle\lnab{}dF(Z,\mu),
JdF(\bmu)\rangle -\langle\lnab{}dF(Z,\bmu),JdF(\mu)\rangle\La{)}.\non
\end{eqnarray}
Thus, since $ Im( iA)= Im (i\bar{A})$ for any complex number $A$,\\[-6mm]
\begin{eqnarray}
\lefteqn{-4\sum_{\be}Im\LA{(}d\kappa_{p_{0}}
\la{(}\lnab{\bbe}\tilde{J}(\be)\la{)}\LA{)}=}\non\\[-3mm]
&\!\!\!\!\!=&\!\!\!\!\!\!
\sum_{\be,\mu}\!\!\frac{32}{\sin^{2}\theta_{\mu}}\LA{(}Im\La{(}
g(\lnab{}dF(\frac{i}{2} \lnab{\bbe}\tilde{J}(\be)
,\bmu),JdF(\mu)\La{)}
-Im\La{(}g\la{(}\lnab{}dF(\frac{i}{2}\lnab{\bbe}
\tilde{J}(\be),\mu),JdF(\bmu)\la{)}\La{)}\LA{)}\non\\[-3mm]
&\!\!\!\!\!=&\!\!\!\!\!\!
\sum_{\be,\mu}\!\!\frac{32}{\sin^{2}\theta_{\mu}}\LA{(}\!\!Im\!\La{(}\!
g(\lnab{}dF(\frac{i}{2} \lnab{\be}\!\tilde{J}(\bbe)
,\mu),JdF(\bmu)\!\La{)}\!-\!Im\!\La{(}\!
g\la{(}\lnab{}dF(\frac{i}{2}\lnab{\bbe}
\!\tilde{J}(\be),\mu),JdF(\bmu)\la{)}\!\La{)}\!\LA{)}
\end{eqnarray}
Consequently, for $F$ minimal, at $p_{0}$,\\[-7mm]
\begin{eqnarray}
\triangle \kappa
&=&4i\sum_{\be} Ricci^{N}(JdF(\be),dF(\bbe))\non\\[-2mm]
&&-\sum_{\be,\mu}\frac{32}{\sin^{2}\theta_{\mu}}Im\La{(} g\la{(}\mbox{\Large$\nabla$}_{\be}(\lnab{}dF)^{(1,1)}(\bbe,\mu),
JdF(\bmu)\la{)}\La{)}\\[-2mm]
&&+\sum_{\be,\mu}\frac{32}{\sin^{2}\theta_{\mu}}Im\La{(} g\la{(}\mbox{\Large$\nabla$}_{\bbe}(\lnab{}dF)^{(1,1)}(\be,\mu),
JdF(\bmu)\la{)}\La{)}\\[-2mm]
&& -4 Im\La{(}d\kappa_{p_{0}}\la{(}\sum_{\be}
\lnab{\bbe}\tilde{J}(\be)\la{)}\La{)}\\[-1mm]
&&-\sum_{\be,\mu,\rho}\!\!\!\frac{64(\cos\theta_{\mu}
\!+\!\cos\theta_{\rho})}{\sin^{2}\theta_{\mu}\sin^{2}\theta_{\rho}}
Re\la{(}\gf{\be}{\mu}{\brho}\gf{\bbe}{\rho}{\bmu}\!\la{)}
\end{eqnarray}
Of course we may assume $k\geq 1$, since, on $\Omega^{0}_{0}$, the
summation on $\mu$ (and $\be$) of the  left-hand side of (2.13) vanish, (because of a symmetry and skew-symmetry argument)
and Proposition 2.1 gives Proposition 2.2 on that open set.
If $F$ is broadly-pluriminimal 
then $F$ is pluriharmonic with respect to any local complex
structure $\tilde{J}$ we described above, on each $\Omega_{2k}^0$.
 Then  $(\lnab{}dF)^{(1,1)}=0$ and (2.22), (2.23),
 and (2.25) vanish. Now we prove that (2.24)
also vanish.  On a neighbourhood of $p_{0}$, since $\tg_{\mu\bmu}=0$, 
for all $\mu\geq k$, then the $\sum_{\mu}$  in (2.20)
and in (2.21) can be replaced by $\sum_{1\leq\mu\leq k}$.
On the other hand, from (2.21), (2.3) and broadly-pluriminimality
\begin{equation}
-4\sum_{\be}Im\LA{(}d\kappa_{p_{0}}
\la{(}\lnab{\bbe}\tilde{J}(\be)\la{)}\LA{)}=\sum_{1\leq \be\leq n}~\sum_{
1\leq \mu\leq k}
\frac{64}{\sin^{2}\theta_{\mu}} Im \La{(}\sum_{1\leq \rho\leq n}
\langle\lnab{\be}\bbe, \brho\rangle \gf{\rho}{\mu}{\bmu}\La{)}
\end{equation}
From broadly-pluriminimality,  for each $\mu\leq k$, $\gf{\brho}{\mu}{\bmu}=0$, $\forall \rho\geq k+1$.
But on ${\cal K}_{\omega}$ 
the $g_{M}$-orthogonal complex structure $ J'$ is arbitrary.
In particular we may replace $J'$ by $-J'$.
This means that we may replace  $\brho\geq k+1$ by $\rho\geq k+1$.
Therefore, $\gf{\rho}{\mu}{\bmu}=0$, 
$\forall \rho\geq k+1$ and $\forall \mu\leq k$.
From Lemma 2.1 and  and pluriharmonicity, for $\rho\leq k$ and $\forall \be$, 
\[ (\cos\theta_{\be}+\cos\theta_{\rho})\langle\lnab{\bbe}{\be},\rho\rangle
=i\gf{\bbe}{\be}{\rho}-i\gf{\bbe}{\rho}{\be}=0\]
with $(\cos\theta_{\be}+\cos\theta_{\rho})>0$. Then $\langle\lnab{\bbe}{\be},\rho\rangle=0$, $\forall \rho\leq k$.
Thus, for each $\mu\leq k$,
\[\langle\lnab{\be}\bbe,\brho\rangle
\gf{\rho}{\mu}{\bmu}=0\]
for any $\rho\leq k$ and also for any $\rho\geq k+1$. That is
$(2.26)=0$, and we have proved Proposition 2.2.~~~~~~~\qed
\begin{Cr} In the conditions of Proposition 2.2, if $N$ is 
K\"{a}hler-Einstein with $ Ricci^{N}= Rg$, on each
$\Omega^{0}_{2k}$, 
\[\triangle \kappa = -2R\La{(}\sum_{\be}\cos\theta_{\be}\La{)}.\] 
\end{Cr}
 \em Remark. \em 
One may search when $(\Omega_{2n},J_{\omega},g_{M})$ is a K\"{a}hler
manifold, for,  it is more usual to define plurharmonicity on K\"{a}hler
manifolds then in almost complex manifolds. More generally, 
we may ask when or where $\lnab{}J_{\omega}=0$
holds. Set, on each $\Omega_{2k}^{0}$,  $\epsilon'_{\al}=+1$, $\epsilon'_{\bal}=-1$, $\forall 1\leq \al \leq k$, and
$\epsilon'_{\al}=\epsilon'_{\bal}=0$, $\forall \al\geq k+1$.
Then,  $~\forall
A,B\in\{1,\ldots,n,\bar{1},\ldots,\bar{n}\},~ \forall Z\in T_{p}M$, and $\forall p\mbox{~near~}p_{0}\in\Omega_{2k}^{0}$\\[-4mm]
\begin{equation}
\langle\lnab{Z}J_{\omega}({A}),{B}\rangle=
i(\epsilon'_{A}+\epsilon'_{B})\langle\lnab{Z}{A},{B}\rangle=-\langle
{A},\lnab{Z}J_{\omega}({B})\rangle.\end{equation}\\[-8mm] 
Then, on $\Omega_{2k}^{0}$\\[-5mm]
\begin{equation}\begin{array}{l}
\lnab{Z}J_{\omega}({\cal K}_{\omega})\subseteq {\cal K}_{\omega}^{\bot}\\
\lnab{Z}J_{\omega}(({\cal K}_{\omega}^{\bot})^{1,0})\!\subseteq\! ({\cal
K}_{\omega}^{\bot})^{0,1}\cup {\cal K}_{\omega}\\ \lnab{Z}J_{\omega}(({\cal
K}_{\omega}^{\bot})^{0,1})\!\subseteq\! ({\cal
K}_{\omega}^{\bot})^{1,0}\cup {\cal K}_{\omega}.\end{array} \end{equation}
For simplicity, we denote by $\cos \theta_{\bbe}=-\cos\theta_{\be}$. Hence
at $p_{0}$ and $\forall A,B\in\{1,\ldots,n,$ $\bar{1},\ldots,\bar{n}\}$, 
$F^{*}\omega({A},{B})=\delta_{A\bar{B}}\frac{i}{2}\cos\theta_{A}$.
Derivating at $p_{0}$
\begin{equation}
d(F^{*}\omega({A},{B}))(Z)=
\lnab{Z}F^{*}\omega({A},{B})-i(\cos\theta_{A}+\cos\theta_{B})
\langle\lnab{Z}{A},{B}\rangle.
\end{equation}
Consequently, from  (2.27) and (2.29), we obtain, at $p_{0}$
\begin{equation}{\small \begin{array}{lcll}
\lnab{Z}F^{*}\omega(\!{A},\!{B}\!)\!&\!\!\! =\!\!\! & 0 & \forall A,
B\in\{k+1,\overline{k+1},\ldots,n,\bar{n}\}\\
\lnab{Z}F^{*}\omega(\!{A},\!{B}\!)\!&\!\!\! =\!\!\! &
\frac{(\cos\theta_{A}+\cos\theta_{B})}
{\epsilon'_{A}+\epsilon'_{B}}\langle\lnab{Z}J_{\omega}({A}),{B}\rangle &
\mbox{whenever~}\epsilon'_{A}+\epsilon'_{B}\neq 0 \end{array}}
\end{equation}
In particular, (2.28), (2.30) and (2.6) let us to conclude that,
 $(\Omega_{2n},J_{\omega},g_{M})$ is a K\"{a}hler manifold
~~iff~~$\lnab{Z}F^{*}\omega$ is of type $(1,1)$, $\forall Z\in T^{c}M$, on
$\Omega_{2n}$ ~~iff~~
$(V,W)\ra g(\lnab{Z}dF(V),JdF(W))$ is symmetric on $T^{1,0}M$,
$\forall Z\in T^{c}M$ on $\Omega_{2n}$
~~iff~~ $\langle\lnab{Z}J_{\omega}({\al}),{\be}\rangle=0$ $~\forall 1\leq
\al,\be\leq n$ and $ ~\forall Z\in T^{c}M$, on $\Omega_{2n}$. 
In [S-V] we prove that,  if $\Fw$ is parallel, that is,
 $ \lnab{Z}F^{*}\omega=0$, then  the K\"{a}lher angles are constant, ${\cal K}_{\omega}$, $ {\cal K}_{\omega}^{\bot}$  are parallel  sub-bundles of $TM$,
and  $\lnab{}J_{\omega}=0$.
\subsection{Some conclusions}
Let $p_{0}$ be an {\bf absolute maximum of} $\kappa$. Then
$\kappa(p_{0})\geq 0$ with equality to zero iff $F^{*}\omega=0$ everywhere,
that is, $F$ is a Lagrangian submanifold. We only know that $\kappa$ is
locally Lipschitz on  $M$ and smooth on each $\Omega^{0}_{2k}$.
Nevertheless we have the following proposition:
 \begin{Pp}
The map $\kappa$ is differentiable at $p_{0}$ with $d\kappa (p_{0})=0$.
\end{Pp}
\em Proof. \em Let $2k$ be the rank of $F^{*}\omega$ at $p_{0}$, and
$\pm i\cos\theta_{\al}(p)$ the eigenvalues of $(F^{*}\omega)^c$,
the $\Co$-linear  extension of $\Fw$ to $T^c M$, for $p$ near
$p_{0}$. Of course we may assume $k\geq 1$. Set 
\[ \kappa_{1}(p)= \sum_{1\leq \al\leq k}\log
\LA{(}\frac{1+\cos\theta_{\al}(p)}{1-\cos\theta_{\al}(p)}\LA{)},~~~~~~~~
\kappa_{2}(p)= \sum_{k+1\leq \al\leq n}\log \LA{(}\frac{1+\cos\theta_{\al}(p)}
{1-\cos\theta_{\al}(p)}\LA{)}.\]
$\kappa_{1}$ is the piece of $\kappa$ defined by angles which cosine is not
zero near $p_{0}$. The re\-maining angles, forming the 
$\kappa_{2}$, are zero
at $p_{0}$, therefore they remain well distinct from the ones which form
the $\kappa_{1}$, for $p$ near $p_{0}$. 
Then we may conclude that $\kappa_{1}$ is
smooth. In fact, if we take $C$ a contour in the 1-complex space around
the eigenvalues $\pm i\cos\theta_{1},\ldots,\pm i\cos\theta_{k}$ for $p$ near $p_{0}$, not meeting any of all other
eigenvalues, then $P=\frac{1}{2\pi i}\int_{C}(\lambda I
-(F^{*}\omega)^{c})^{-1}d\lambda$, where $I:T^c M\ra T^c M$ is the
identity morphism, is, for each $p$ on a neighbourhood of $p_{0}$,
 an orthogonal projection (since $F^{*}\omega$ is normal)
onto the subspace that is the direct sum of the eigenspaces corresponding
to the eigenvalues inside $C$, and $T =\frac{1}{2\pi i}\int_{C}\lambda
(\lambda I -(F^{*}\omega)^{c})^{-1}d\lambda$ is the restriction of
$(F^{*}\omega)^{c}$ to the sub-bundle $E= P(T^{c}M)$. Clearly $T$ and $P$
are smooth sections of $T^{*}M^{c}\otimes T^{c}M$, and so $E$ is smooth.
Then $\kappa_{1}$ is
just defined as $\kappa$, but relative to $T$ on $E$ where it is 
invertible everywhere. Now we have, for $p$ near $p_{0}$, 
$\kappa=\kappa_{1}+\kappa_{2}$,  with $\kappa_{1}$ smooth.
The point $p_{0}$ is a maximum of $\kappa$ and of $\kappa_{1}$, and is a
minimum of $\kappa_{2}$. Then $d\kappa_{1}(p_{0})=0$. We may assume
$\kappa$ and $\kappa_{i}$ are defined on a open set of $\R{m}$. We will
prove that $\kappa_{2}$ is stationary at $p_{0}$, that is, $\kappa_{2}$ is
differentiable at $p_{0}$ and its differential at $p_{0}$ is zero. Thus we
want to prove that $r(p)= \kappa_{2}(p)-\kappa_{2}(p_{0})$ satisfies
\\[-5mm] 
\begin{equation}
 \forall t>0 \exists s>0: ||p-p_{0}||\leq s \Longrightarrow
||r(p)||\leq t ||p-p_{0}||
\end{equation}
\\[-6mm]
 Suppose this is not true.
Then \\[-5mm]
\begin{equation}
\exists t'>0\mbox{~and~}p_{n}\ra p_{0}\mbox{~~such~that~}
||r(p_{n})||>t'||p_{n}-p_{0}||
\end{equation}
\\[-6mm] 
From $\kappa(p_{n})\leq
\kappa(p_0)$,
$\kappa_{1}(p_{n})\leq \kappa_{1}(p_0)$, and
$\kappa_{2}(p_{n})\geq\kappa_{2}(p_0)$, we have 
\begin{eqnarray*}
0&\leq& \frac{r(p_{n})}{||p_{n}-p_{0}||}
=\frac{\kappa_{2}(p_{n})-\kappa_{2}(p_0)}{||p_{n}-p_{0}||}
=\frac{\kappa(p_{n})-\kappa(p_0)}{||p_{n}-p_{0}||}
-\frac{\kappa_{1}(p_{n})-\kappa_{1}(p_0)}{||p_{n}-p_{0}||}\\
&\leq&-\frac{\kappa_{1}(p_n)-\kappa_{1}(p_0)}{||p_{n}-p_{0}||}=
-d\kappa_{1}(p_{0})\LA{(}\frac{p_{n}-p_{0}}{||p_{n}-p_{0}||}\LA{)}
+\frac{r'(p_{n})}{||p_{n}-p_{0}||}
\end{eqnarray*}
where $r'(p)$ satisfies (3.31). Since $d\kappa_{1}(p_{0}) =0$, the last term
in the above inequality converges to zero, contradicting
$(3.32)$.~~~~~~\qed\\[1mm]

In the conditions of Proposition 2.2 and Corollary 2.1 assume that $M$ is compact and $N$ is a K\"{a}hler-Einstein manifold. If $R<0$, 
by applying maximum principle to $\kappa$ at a maximum point $p_{0}$, we immediately conclude from Corollary 2.1:
\begin{Lm} If $R<0$, and $F$ is broadly-pluriminimal, not Lagrangian and 
with no complex directions, then $p_{0}$ is not in $\Omega_{2k}^{0}$
$\forall 0\leq k\leq n$. In other words, if the rank of $F^{*}\omega$ is 2k 
at $p_{0}$, then there exists a sequence $p_{m}\ra p_{0}$ such that the rank
of $F^{*}\omega$ at $p_{m}$ is $>2k$. In particular $F^{*}\omega$ cannot
have constant rank and
$F^{*}\omega$ is degenerated at $p_{0}$. \end{Lm}
\begin{Pp} If $M$ is compact and  $R<0$, 
every broadly-pluriminimal immersion either
has Lagrangian or complex directions.
\end{Pp}
\em Proof. \em If $F$ has no Lagrangian directions
$M=\Omega_{2n}=\Omega^{0}_{2n}$. Then, non existence of complex directions
contradicts the Lemma.~~~~~~\qed \\[6mm] 
 Of course if $n=1$ the assumption $(ii)$ for broadly-pluriminimal is automatic,  and minimal surfaces are broadly-pluriminimal, and we recover some result of Wolfson [W].
\begin{Pp} If $F:M\ra
N$ is broadly-pluriminimal  without complex directions, 
$M$ is compact, $R<0$ and if rank $F^{*}\omega\leq 2$, 
then $F$ is Lagrangian.
\end{Pp}
\em Proof. \em If $\kappa(p_{0})\neq 0$, then $F^{*}\omega$ would have rank
2 at $p_{0}$, contradicting Lemma 2.2.\qed \\[3mm]
An immediate consequence of Corollary 2.1 by making $\kappa$ constant, 
is the following proposition:
\begin{Pp} If $F$ is broadly-pluriminimal, not Lagrangian, with constant K\"{a}hler angles and no complex directions, then $R=0$ must hold.
\end{Pp}
The following proposition is  anounced  in [S-V], but only proved in
this paper, as a consequence of Corollary 2.1 as well:
\begin{Pp} If $M$ is compact, $N$ is K\"{a}hler-Einstein with $R<0$, and if
F is broadly-pluriminimal  without complex directions and
with equal K\"{a}hler angles, then F is Lagrangian. \end{Pp}
\em Proof. \em If the K\"{a}hler angles are all equal, $\cos\theta_{\al}=
\cos\theta, \forall\al$, and if this is not zero everywhere, that is, F is not Lagrangian, then $p_{0}\in \Omega_{2n}$, what is not possible.~~~~\qed 
\begin{Th} If $n=2$, and if $F:M\ra N$ is broadly-pluriminimal 
 without complex directions, $M$ is
compact and orientable, and $N$ is K\"{a}hler-Einstein with $R<0$, then $F$
is Lagrangian. \end{Th}
\em Proof. \em We have two eigenvalues $\cos\theta_{1}\geq
\cos\theta_{2}\geq 0$, and we suppose $F$ not Lagrangian. Then
from Lemma 2.2 we must have $\cos\theta_{1}(p_{0})>0$,
$\cos\theta_{2}(p_{0})=0$ and there exist $p_{m}\ra p_{0}$ such that
$\cos\theta_{2}(p_{m})>0$, $\forall m$. Since $\cos\theta_{1}\neq
\cos\theta_{2}$ at $p_{0}$, by continuity they remain
different for $p$ near $p_{0}$, and
$\pm i \cos\theta_{1}$ do not vanish and have multiplicity one 
for $p$ near
$p_{0}$. Therefore, an application of the implicit map theorem to the characteristic equation of the complex extension of $\Fw$, shows that $\cos\theta_{1}$ defines a smooth map defined on a neighbourhood of
$p_{0}$. From $F^{*}\omega (p)=\cos\theta_{1}(p) X_{*}^{1}\wedge
Y_{*}^{1}+\cos\theta_{2}(p)X_{*}^{2}\wedge Y_{*}^{2}$, where
$X_{1},Y_{1},X_{2},Y_{2}$
is a diagonalizing orthonormal basis of $F^{*}\omega(p)$, we have
\[(F^*\omega)^{2}(p)= \epsilon(p)2\cos\theta_{1}(p)\cos\theta_{2}(p)Vol\]
where $Vol$ is the volume element of $M$, and where $\epsilon (p)$ is $+1$
or $-1$  according
 $X_{1},Y_{1},X_{2},Y_{2}$ is an direct or inverse basis,
respectively. Since $(F^*\omega)^{2}$ and $\cos\theta_{1}$ are smooth, then
\begin{equation}
s_{2}(p)=\epsilon(p)\cos\theta_{2}(p)
\end{equation}
is a smooth map for $p$ near $p_{0}$.
Without loss of generalization we may suppose that $s_{2}(p_{m})>0$
$\forall m$, by taking a subsequence, and  changing the orientation of
$M$, if necessary. 
Then $\epsilon(p_{m})=1$ $\forall m$. Let us take the smooth map,
defined on a neigbourhood of $p_{0}$\\[-4mm] 
\[\tilde{\kappa}=
\log\left(\frac{1+\cos\theta_{1}}{1-\cos\theta_{1}}\right)
+\log\left(\frac{1+s_{2}}{1-s_{2}}\right)\]\\[-4mm]
 Since $s_{2}=\cos\theta_{2}$ in the open set of the points where 
$s_{2}>0$ or
equivalently where $\epsilon =+1$ and $\cos\theta_{2}$ does not vanish,
then $\kappa =\tilde{\kappa}$ on that set, and in particular in a
neigbourhood of each $p_{m}$, for each $m$. Moreover $\kappa(p_{0})=
\tilde{\kappa}(p_{0})$.
We now prove that $p_{0}$
is also a maximum of $\tilde{\kappa}$.
Let $p$ on a neigbourhood of $p_{0}$. If $s_{2}(p)=\cos\theta_{2}(p)$, then
$\tilde{\kappa}(p)=\kappa(p)\leq\kappa(p_{0})=\tilde{\kappa}(p_{0})$. If
$s_{2}(p)=-\cos\theta_{2}(p)$ then\\[-3mm]
\[\tilde{\kappa}(p)=\kappa(p)-2\log\left(\frac{1+\cos\theta_{2}(p)}
{1-\cos\theta_{2}(p)}\right)\]\\[-4mm]
with $\log\left(\frac{1+\cos\theta_{2}(p)} {1-\cos\theta_{2}(p)}\right)\geq 0$.
From $\kappa(p)\leq \kappa(p_{0})=\tilde{\kappa}(p_{0})$ we obtain\\[-2mm] 
\[
\tilde{\kappa}(p)=\kappa(p)-2\log\left(\frac{1+\cos\theta_{2}(p)}
{1-\cos\theta_{2}(p)}\right)\leq
\tilde{\kappa}(p_{0})-2\log\left(\frac{1+\cos\theta_{2}(p)}
{1-\cos\theta_{2}(p)}\right)\leq \tilde{\kappa}(p_{0}).\]\\[-3mm] 
Therefore, by maximum principle applied to $\tilde{\kappa}$ at $p_{0}$,
by Corollary 2.1 applied to $p_{m}\in \Omega_{2n}^{0}$, and 
by  continuity of the maps $\cos\theta_{\al}$,\\[-7mm]
\begin{eqnarray*}
0\geq\!\!&&\!\!\!\!\!\triangle \tilde{\kappa}(p_{0})=\lim_{m}\triangle
\tilde{\kappa}(p_{m})=\lim_{m}\triangle \kappa(p_{m})= \lim_{m} -2
R(\cos\theta_{1}(p_{m})+\cos\theta_{2}(p_{m}))\\ 
&&=-2R(\cos\theta_{1}(p_{0})+\cos\theta_{2}(p_{0}))= 
-2R\cos\theta_{1}(p_{0})\geq 0
\end{eqnarray*}\\[-8mm]
which implies $\cos\theta_{1}(p_{0})=0$, that is a contradiction.
~~~~\qed\\[5mm]
Unfortunately a similar argument does not work for higher dimensions, since
we do not have smooth functions as good as (2.33). 
We only have smoothness of some
symmetric functions of the K\"{a}hler angles. \\[8mm]
{\Large \bf{References}}\\[2mm]
[Bh]~~~~{\small R.\ Bhatia, \em
Matrix Analysis\em, GTM {\bf 169}, Springer, 1996.}\\[0mm]
[S-V]~~{\small I.M.C.\ Salavessa \& G.\ Valli, \em Minimal
submanifolds of K\"{a}hler-Einstein manifolds with equal K\"{a}hler angles,
\em, e-print no. math.DG/0002050 of 7 Feb 2000, submited to publication}\\[0mm]
[O-V]~~{\small Y.\ Ohnita \& G.\ Valli \em Pluriharmonic maps into
 compact Lie groups and factorization into unitons, \em
Proc.\ London Math.\ Soc. {\bf 61} (1990), 546-570.}\\[0mm] 
[W]~~~~~{\small J.\ G.\ Wolfson, \em Minimal
Surfaces in K\"{a}hler Surfaces and Ricci Curvature, \em J. Diff. Geom,
{\bf 29}, 281--294, (1989).}\\[4mm] 
\end{document}